\documentclass[a4paper]{article}
\usepackage[top=25truemm,bottom=25truemm,left=20truemm,right=20truemm]{geometry}
\usepackage{amsmath,amsthm,amsfonts,amssymb,stmaryrd,indentfirst,tabularray,makecell,hhline,comment}
\usepackage[mathscr]{euscript}
\usepackage[shortlabels]{enumitem}
\usepackage[labelfont=rm]{caption}
\usepackage{titlesec}
\usepackage{tikz-cd}
\usepackage{url}
\usepackage[backref]{hyperref}
\hypersetup{
    colorlinks = true,
    linkcolor = [rgb]{0.3,0.55,0.25},
    citecolor = [rgb]{0.164,0.39,0.77},
    urlcolor = [rgb]{0.7,0.2,0.4}
}

\usepackage{tikz}
\usetikzlibrary{decorations.markings}

\newcommand{\RGbar}{
	\begin{tikzpicture}[line width=0.6pt, baseline=-1ex, scale=1.2]
	
	\node[label=center:{\footnotesize$\bar 1$}] at (0.25,-0.25)  {};
	
	\node[fill, circle, inner sep=0pt, minimum size=3pt] at (0,0) (n1) {};
	\node[fill, circle, inner sep=0pt, minimum size=3pt] at (0.5,0) (n2) {};
	
	\node[label=center:{\footnotesize$1$}] at (0,0.2) {};
	\node[label=center:{\footnotesize$2$}] at (0.5,0.2) {};
	
	\draw[-Stealth] (n1)--(n2);
	
	\end{tikzpicture}
}

\newcommand{\RGtadpole}{
	\begin{tikzpicture}[line width=0.6pt, baseline=-1ex, scale=1.2]
	
	\node[label=center:{\footnotesize$\bar 2$}] at (-0.6,0.2)  {};
	\node[label=center:{\footnotesize$\bar 1$}] at (-0.2,0)  {};
	
	\node[fill, circle, inner sep=0pt, minimum size=3pt] at (0,0) (n1) {};
	
	\node[label=center:{\footnotesize$1$}] at (0.2,0) {};
	
	\draw (n1) to[out=90, in=0] (-0.2,0.22) to[out=180, in=90] (-0.4,0) to[out=-90, in=180] (-0.2,-0.2);
	\draw[-Stealth] (-0.2,-0.2) to[out=0, in=238] (n1);
	
	\end{tikzpicture}
}

\newcommand{\RGloop}{
	\begin{tikzpicture}[line width=0.6pt, baseline=-1ex, scale=1.2]
	
	\node[label=center:{\footnotesize$\bar 1$}] at (0,0)  {};
	\node[label=center:{\footnotesize$\bar 2$}] at ({3.7*360/7}:1.4)  {};
	
	\node[fill, circle, inner sep=0pt, minimum size=3pt] at ({0*360/7}:1) (n0) {};
   	\node[fill, circle, inner sep=0pt, minimum size=3pt] at ({1*360/7}:1) (n1) {};
	\node[fill, circle, inner sep=0pt, minimum size=3pt] at ({2*360/7}:1) (n2) {};
   	\node[fill, circle, inner sep=0pt, minimum size=3pt] at ({3*360/7}:1) (n3) {};
	\node[inner sep=0pt,minimum size=12pt] at ({4*360/7}:1) (n4) {};
   	\node[fill, circle, inner sep=0pt, minimum size=3pt] at ({5*360/7}:1) (n5) {};
	\node[fill, circle, inner sep=0pt, minimum size=3pt] at ({6*360/7}:1) (n6) {};
	
	\node[label=center:{\footnotesize$1$}] at ({0*360/7}:1.2) {};
	\node[label=center:{\footnotesize$2$}] at ({1*360/7}:1.2) {};
	\node[label=center:{\footnotesize$3$}] at ({2*360/7}:1.2) {};
	\node[label=center:{\footnotesize$4$}] at ({3*360/7}:1.2) {};
	\node[label=center:{\;\footnotesize$k\!-\!1$}] at ({5*360/7}:1.2) {};
	\node[label=center:{\footnotesize$k$}] at ({6*360/7}:1.2) {};
	
	\draw[-Stealth] (n0)--(n1);
	\draw[-Stealth] (n1)--(n2);
	\draw[-Stealth] (n2)--(n3);
	\draw[dashed] (n3)--(n4);
	\draw[dashed, -Stealth] (n4)--(n5);
	\draw[-Stealth] (n5)--(n6);
	\draw[-Stealth] (n6)--(n0);
	
	\end{tikzpicture}
}

\DeclareMathOperator{\id}{id}

\DeclareMathOperator{\Der}{Der}

\DeclareMathOperator{\Hom}{Hom}

\DeclareMathOperator{\End}{End}

\DeclareMathOperator{\Tr}{Tr}
\DeclareMathOperator{\Div}{Div}

\DeclareMathOperator{\sdiv}{div} % stands for single div

\DeclareMathOperator{\DR}{DR}

\DeclareMathOperator{\ad}{ad}

\newcommand{\pair}[2]{\langle#1,#2\rangle}
\newcommand{\op}[0]{^\mathrm{op}}

\theoremstyle{plain}
\newtheorem{theorem}{Theorem}
\newtheorem*{theorem*}{Theorem}
\newtheorem{proposition}[theorem]{Proposition}
\newtheorem*{proposition*}{Proposition}

\newtheorem{lemma}[theorem]{Lemma}
\newtheorem{corollary}[theorem]{Corollary}

\newtheorem{question}[theorem]{Question}

\theoremstyle{definition}
\newtheorem{definition}[theorem]{Definition}
\newtheorem*{definition*}{Definition}
\newtheorem{remark}[theorem]{Remark}
\newtheorem{example}[theorem]{Example}

\numberwithin{theorem}{section}
\everymath{\displaystyle}
\allowdisplaybreaks

\usepackage{apptools}
\titleformat{\section}{\large\scshape}{\IfAppendix{\appendixname}{} \thesection.}{3pt}{}
\titleformat{\subsection}{\scshape}{\thesubsection.}{3pt}{}

\usepackage{titling}  
\pretitle{
	\begin{center}\Large\bfseries
}
\posttitle{
	\end{center}\ \\[-20pt]
}
\preauthor{
	\begin{center}\large
}
\postauthor{
	\end{center}\ \\[-60pt]
}

\title{A Family of Algebraic Operations\\ Extending the Turaev Cobracket}
\author{Toyo TANIGUCHI \thanks{Graduate School of Mathematical Sciences, The University of Tokyo. 3-8-1, Komaba, Meguro-ku, Tokyo, 153-8914, Japan. E-mail: \texttt{toyo(at)ms.u-tokyo.ac.jp}}}
\date{}

\begin{document}
\maketitle

\begin{abstract}
\noindent We introduce a family of maps parametrised by certain ribbon graphs. It is based on a connection in non-commutative geometry and contains the double divergence as a special case. Applying the construction to the case of the group algebra of the fundamental group of a compact connected oriented surface with boundary, we obtain an algebraic generalisation of the Turaev cobracket. If the connection is flat, they define classes in the Lie algebra cohomology of the space of derivations. In the case of the free associative algebra, we show that they are canonically identified with the standard generators of the cohomology ring of the matrix Lie algebra $\mathfrak{gl}_n$.\\
\end{abstract}

\noindent{\textit{2020 Mathematics Subject Classification: 16D20, 17B56, 17B65, 53C05, 57K20, 58B34.}}\\
\noindent{\textbf{Keywords:} non-commutative geometry, divergence maps, flat connections, ribbon graph, loop operations.}

\section{Introduction}

The Turaev cobracket is a \textit{loop operation} introduced in \cite{turaev} together with its quantisation, and gives a Lie bialgebra structure on the vector space spanned by the homotopy set of free loops on a connected compact oriented surface. This Lie bialgebra is a topological counterpart of the Batalin-Vilkovisky structure on the moduli space of flat principal super Lie group bundles; the Lie algebra part is classically known by Goldman \cite{goldman}, and the extension to the whole BV structure is done by Alekseev, Naef, Pulmann, and \v{S}evera \cite{anps}.

If the surface has only one boundary component, the associated graded of this Lie bialgebra with a suitable filtration is isomorphic to Schedler's necklace Lie bialgebra on the trace space $|T(W)|$ of the free associative algebra over a finite-dimensional vector space $W$; see Section 6 of \cite{gtjohnson} for the precise statement and original references. This algebraic operation also appears in the framework of operad theory. In \cite{ribbon}, Merkulov--Willwacher introduced a prop of ribbon graphs with a canonical representation on $|T(W)|$ where $W$ is equipped with a skew-symmetric pairing. In particular, Theorem 4.2.3 in \cite{ribbon} shows that a representation over this prop is automatically an involutive Lie bialgebra, and it can be shown that their bialgebra structure on $|T(W)|$ coincides with Schedler's one.\\

In this paper, we will define a family of algebraic operations on the space of free loops extending the Turaev cobracket. We discuss the construction in a special case. First, we define a series of maps called $k$-\textit{divergence}
\[
	\Div^{\nabla}_k\colon {\Der_\mathbb{K}(A)}^{\otimes k} \to |A|^{\otimes 2}
\]
associated with a $\mathbb{K}$-algebra $A$, an $A$-bimodule $M$, an action of $\Der_\mathbb{K}(A)$ on $M$ and a non-commutative version of a connection (or a covariant derivative) $\nabla$ on $M$. In some situations, $|A|$ is endowed with a  Lie bracket and a Lie algebra homomorphism $\psi \colon |A| \to \Der_\mathbb{K}(A)$. The restriction
\[
	\Div^\nabla_k\circ \,\psi^{\otimes k}\colon |A|^{\otimes k} \to |A|^{\otimes 2}
\]
is denoted by $\delta^{\psi,\nabla}_k$.

For $k=1$, the map $\Div^{\nabla}_1$ is an analogue of the usual divergence; see the author's previous paper \cite{toyo} for more details. Its pull-back $\delta^{\psi,\nabla}_1$ has a very interesting interaction with loop operations. The result of Alekseev--Kawazumi--Kuno--Naef \cite{akkn} and the author \cite{toyo} state that, in the case of group algebra $A = \mathbb{K}\pi$ of the fundamental group $\pi = \pi_1(\Sigma)$ of a connected oriented compact surface $\Sigma$ with non-empty boundary, we can choose $M$, $\nabla$ and $\psi$ suitably so that $-\delta^{\psi,\nabla}_1$ coincides with the framed version of the Turaev cobracket in \cite{akkn}. In fact, the Lie algebra map $\psi$ is given by the Kawazumi--Kuno action $\sigma$ introduced in \cite{sigma}, which is another loop operation. See Remark \ref{rem:cob} for a more precise statement. For $k=2$, some values of $\delta_2^{\sigma,\nabla}$ are computed in Example \ref{ex:second}.

The original motivation for the construction of the $k$-divergences was to craft another loop operation such that multiple intersections are resolved at the same time. Our rather algebraic definition of the operations $\delta^{\psi, \nabla}_k$ lacks topological perspective: even in the case of the group algebra $\mathbb{K}\pi$, each operation is not mapping-class-group equivariant as opposed to the case of the Turaev cobracket, but instead, the framework is applicable to an arbitrary algebra. See Remark \ref{rem:unsuc} for the discussion.\\

In the case of $A = T(W)$, the free associative algebra over a finite-dimensional $\mathbb{K}$-vector space $W$, the ribbon graph operations reappear. Suppose we are given a \textit{skew-symmetric} pairing
\[
	\pair{\cdot}{\cdot}\colon W\otimes W \to \mathbb{K},
\]
from which we obtain a Lie algebra homomorphism $\mathrm{Ham}_{\pair{\cdot}{\cdot}}\colon |T(W)| \to \Der_\mathbb{K}(T(W))$. Here, $|T(W)|$ is equipped with the necklace bracket. 

\begin{theorem*}[Theorem \ref{thm:ribbonop}]
Let $\nabla\!_W$ be the canonical flat connection given in Definition \ref{def:connW} on the space of non-commutative $1$-forms on $T(W)$. Then, the map
\[
	(-1)^k\delta^{\mathrm{Ham}_{\pair{\cdot}{\cdot}},\nabla\!_W}_k\colon |T(W)|^{\otimes k} \to |T(W)|^{\otimes 2}
\]
coincides with the operation associated with the unique bivalent connected ribbon graph with $k$ vertices and $k$ edges.
\end{theorem*}

Recall that the divergence map on an oriented Riemannian manifold in the usual differential geometry is a $1$-cocycle on the Lie algebra of vector fields. The following is a corresponding result for $k$-divergences. Let $A$ be any $\mathbb{K}$-algebra and $\mathrm{alt}\colon \wedge^k \Der_\mathbb{K}(A) \to \Der_\mathbb{K}(A)^{\otimes k}$ the anti-symmetrisation map.

\begin{theorem*}[Theorem \ref{thm:kcocyc}]
If the connection $\nabla$ is flat, the map $\Div^\nabla_k\!\circ\,\mathrm{alt}$ is a Lie algebra $k$-cocycle on $\Der_\mathbb{K}(A)$. If $k$ is even, the cocycle vanishes.
\end{theorem*}
\noindent Along the way, we develop a theory of non-commutative de Rham cohomology with vector bundle coefficients, which is eventually applied to give a proof of the theorem above.\\

In the important case of $A = T(W)$, we have a non-vanishing of the cocycle for $k$ odd. 
\begin{theorem*}[Theorem \ref{thm:dercohom}]
The cohomology class of $\Tr(c_{\nabla\!_W}^{\;\;\;\,k})$ does not vanish $H_\mathrm{CE}^k(\Der_\mathbb{K}(T(W)),|T(W)^\mathrm{e}|)$ for odd $k$.
\end{theorem*}

\noindent\textbf{Organisation of the paper.} In Section \ref{sec:higher}, We define the $k$-divergences and calculate some values in the case of a group algebra. Section \ref{sec:ribbonop} gives the proof of Theorem \ref{thm:ribbonop}. Section \ref{sec:tdr} is devoted to the general theory of non-commutative de Rham cohomology with vector bundle coefficients, which is used in the proof of Theorem \ref{thm:kcocyc}. In Appendix \ref{sec:infact}, we define the infinitesimal action of a derivation to a connection and discuss the relation with the associated divergence for further analysis of divergence maps.\\

\noindent\textbf{Acknowledgements.} The author thanks Sergei Merkulov for the heart-warming discussion and for sharing their draft on ribbon graphs and loop operations, Nariya Kawazumi for many comments on the draft and Geoffrey Powell for pointing out a critical error in the previous version.\\

\noindent\textbf{Conventions.} $\mathbb{K}$ is a field of characteristic zero. Unadorned tensor products are always over $\mathbb{K}$. We use Sweedler-type notation: an element of the tensor product $X\otimes X$ is denoted as $x = x'\otimes x''$.\\

\section{Divergence Maps and a Family of Operations}\label{sec:higher}

Let $B$ be a unital associative $\mathbb{K}$-algebra and $B^\mathrm{e} = B\otimes B\op$ the universal enveloping algebra. Recall some notions introduced in \cite{toyo}. Also, see \cite{ginz} for a general theory of non-commutative differential geometry.

\begin{definition}\ 
\begin{itemize}
	\item The \textit{trace space} of $B$ is defined by $|B| = B/[B, B]$, which is the quotient by the $\mathbb{K}$-vector space generated by commutators. This is merely the cyclic quotient, not the abelianisation.
	\item The space of \textit{non-commutative $1$-forms} $\Omega^1B$ on $B$ is the $B$-bimodule defined by
	\[
		\Omega^1B = B^\mathrm{e}\{db\colon b\in B\}/\langle d(bb') = (db)b' + b(db'), d\textrm{ is }\mathbb{K}\textrm{-linear}\rangle.
	\]
	Note that every element can be written in the form $bdb'$ by the Leibniz rule. The space of $p$-forms is defined by $\Omega^pB =  \underbrace{\Omega^1B\otimes_B\cdots\otimes_B\Omega^1B}_{p\textrm{ times}}$.
	\item For $f\in\Der_\mathbb{K}(B)$, the \textit{Lie derivative} $L_f\colon \Omega^1B \to \Omega^1B$ is defined by 
	\[
		L_f(bdb') = f(b)db' + bdf(b')\quad \textrm{for } b,b'\in B.
	\]
	\item A \textit{derivation action} on a $B$-module $M$ is a triple $(\mathfrak d, \varphi, \rho)$ where $\mathfrak d$ is a Lie algebra over $\mathbb{K}$ and
	\[
		\varphi\colon \mathfrak d \to \Der_\mathbb{K}(B)\;\textrm{ and } \;\rho\colon\mathfrak d \to \End_\mathbb{K}(M)
	\]
	are Lie algebra homomorphisms satisfying the compatibility:
	\[
		\rho(f)(bm) = \varphi(f)(b) m + b\rho(f)(m)\quad\textrm{for } f\in\mathfrak d, \,b\in B\textrm{ and } m\in M\,.
	\]
	\item A \textit{connection} on a $B$-module $M$ is a $\mathbb{K}$-linear map $\nabla\colon M\to \Omega^1B\otimes_BM$ satisfying the Leibniz rule:
	\[
		\nabla(bm) = db\otimes m + b\nabla(m)\quad\textrm{for } b\in B \textrm{ and }m\in M.
	\]
	This is extended to the degree $+1$ map by the following:
	\begin{align*}
		\nabla\colon \Omega^pB\otimes_BM& \to \Omega^{p+1}B\otimes_BM\\
		\omega\otimes m \hspace{9.5pt}&\mapsto d\omega\otimes m + (-1)^p \omega\nabla(m),
	\end{align*}
	and the connection $\nabla$ is \textit{flat} if $\nabla^2\colon M \to \Omega^2B\otimes_BM$ vanishes.
	\item A $B$-module $M$ is said to be \textit{dualisable} if it is finitely geneated and projective. This is equivalent to the following: for any $B$-module $N$, the natural map $M^*\otimes_{B}N\to \Hom_{B}(M,N)$ is an isomorphism. Here $N^*= \Hom_B(N,B)$ is the $B$-dual of $N$, which is a right $B$-module.
	\item For a dualisable $B$-module $M$, the trace map is defined by the composite
	\begin{align*}
		\Tr\colon \End_B(M) \cong M^*\otimes_BM &\to |B|\\
		\mu\otimes m\hspace{8.5pt}&\mapsto |\mu(m)|.
	\end{align*}
	\end{itemize}
\end{definition}

\begin{definition}
Let $(\mathfrak d, \varphi, \rho)$ be a derivation action on $M$, and $\nabla$ a connection on $M$. With these data given, we define a $\mathbb{K}$-linear map $c_\nabla$ as
\begin{align*}
	c_\nabla\colon \mathfrak d& \to \End_B(M)\\
	 f &\mapsto (i_{\varphi(f)}\otimes \id_M)\circ\nabla - \rho(f)\,,
\end{align*}
where $i_{\varphi(f)}\colon \Omega^1B \to B$ is the contraction, which is a $B^\mathrm{e}$-module map. Note that it takes its value in the space of $B$-linear maps by the compatibility.
\end{definition}

The properties of this map are discussed in Section \ref{sec:tdr}.

\begin{definition}\label{def:divk}
Suppose $M$ is dualisable (i.e., finitely generated and projective). For $k\geq 1$, we define the $k$-\textit{divergence map} $\Div^{\nabla}_k\colon {\mathfrak d}^{\otimes k} \to |B|$ associated with a derivation action $(\mathfrak d, \varphi, \rho)$ and a connection $\nabla$ by the following: for $f_1,\dotsc,f_k\in \mathfrak d$, 
\begin{align*}
	\Div^{\nabla}_k(f_1,\dotsc,f_k) = \Tr(c_\nabla(f_1)\circ\cdots\circ c_\nabla(f_k))\,.
\end{align*}
\noindent We denote the restriction by a Lie algebra homomorphism $\psi\colon \mathfrak{g} \to \mathfrak d$ for some Lie algebra $\mathfrak g$ by $\delta^{\psi,\nabla}_k = \Div^\nabla_k\circ \,\psi^{\otimes k}$.
\end{definition}

\begin{proposition}
$\Div_k^\nabla$ is cyclically symmetric in its arguments: $\Div^\nabla_k(f_1,f_2,\dotsc,f_k) = \Div^\nabla_k(f_2,\dotsc,f_k,f_1)$.
\end{proposition}
\noindent Proof. This follows from the property of the trace map.\qed

\begin{remark}
In the case of $k=1$, this is simply called the \textit{divergence map}, and its details and applications are discussed in \cite{toyo}. Note that the divergence here differs from \cite{toyo} by a sign; this is just a choice of convention.
\end{remark}

\begin{definition}
For a $\mathbb{K}$-algebra $A$, the default derivation action is defined by $B = A^\mathrm{e}$, $M = \Omega^1A$, $\mathfrak d = \Der_\mathbb{K}(A)$, 
\begin{align*}
	\varphi&\colon \Der_\mathbb{K}(A) \to \Der_\mathbb{K}(A^\mathrm{e})\hspace{9.5pt}: f\mapsto f\otimes\id + \id\otimes f\,,\textrm{ and}\\
	\rho&\colon \Der_\mathbb{K}(A)\to\End_\mathbb{K}(\Omega^1A):f \mapsto L_f\,,
\end{align*}
and we will use this action unless otherwise specified.
\end{definition}

\begin{remark}\label{rem:cob}
Let $\Sigma = \Sigma_{g,n+1}$, an oriented connected compact surface of genus $g$ with $n+1$ boundary components. Put $\pi = \pi_1(\Sigma,*)$, the fundamental group of $\Sigma$ with a base point $\ast\in\partial\Sigma$, which is a free group, and consider the group algebra $A = \mathbb{K}\pi$. Then, the Kawazumi-Kuno action $\sigma\colon |\mathbb{K}\pi| \to \Der_\mathbb{K}(\mathbb{K}\pi)$ is defined by the following formula: for a free loop $\alpha$ and a based loop $\beta$, both generically immersed, we set
\[
	\sigma(\alpha)(\beta) = \sum_{p\in\alpha\cap\beta} \mathrm{sign}(p;\alpha,\beta) \,\alpha*_p\beta
\]
where $\mathrm{sign}(p;\alpha,\beta)$ is a local intersection number with respect to the orientation of $\Sigma$, and $\alpha *_p\beta$ is a curve obtained by surgery at the intersection $p$. This is a Lie algebra homomorphism when $|\mathbb{K}\pi|$ is endowed with the structure of the Goldman Lie algebra. See \cite{akkn} for the details.

Now take a free generating system $\mathcal{C} = (\alpha_i,\beta_i,\gamma_j)_{1\leq i\leq g,1\leq j\leq n}$ of $\pi$ so that $\alpha_i$ and $\beta_i$ form a genus pair, $\gamma_j$ is a boundary loop and
\[
	\alpha_1\beta_1\alpha_1^{-1}\beta_1^{-1}\cdots \alpha_g\beta_g\alpha_g^{-1}\beta_g^{-1}\gamma_1\cdots\gamma_n
\]
represents a boundary loop based at $*\in\partial \Sigma$. The flat connection $\nabla\!_\mathcal{C}$ on $\Omega^1\mathbb{K}\pi$ is defined by $\nabla\!_\mathcal{C}((dc)c^{-1}) = 0$ for all $c\in \mathcal{C}$; see author's previous paper \cite{toyo} for the detail. Then, by the result of Alekseev--Kawazumi--Kuno--Naef \cite{akkn} and the author \cite{toyo}, the map $-\delta^{\sigma,\nabla\!_\mathcal{C}}_1\colon |\mathbb{K}\pi|\to |\mathbb{K}\pi|^{\otimes 2}$ coincides with the framed version of the Turaev cobracket. From this perspective, the series of maps $\delta^{\sigma,\nabla\!_\mathcal{C}}_k$ is regarded as an algebraic extension of the Turaev cobracket.
\end{remark}

\begin{table}[htb]
\SetTblrInner{rowsep=1.5pt}
\centerline{\begin{tblr}{|c|c||l|c|}\hline
$x$ & $y$ & $\delta^{\sigma,\nabla\!_\mathcal{C}}_2(x,y)$ & $(x,y)$ \!disjoint?\\\hline
$s^m$ & any element & $0$ & Y\\\hline
$st$ & $st$ & $3stst\otimes 1- 2tsst\otimes 1 - 2st\otimes st + 1\otimes stst$ & N\\\hline
$st$ & $tu$ & $tus\otimes t - tuts\otimes 1 + tust\otimes 1 - tu\otimes st - st\otimes tu + t\otimes t^{-1}sttu$ & N\\\hline
$st$ & $tu^{-1}$ & \makecell{\hspace{-46pt}$tu^{-1}s\otimes t - tu^{-1}ts\otimes 1+ tu^{-1}st\otimes 1$\\\hspace{30pt}$- tu^{-1}\otimes st - st\otimes tu^{-1} + t\otimes t^{-1}sttu^{-1}$} & Y\\\hline
$su$ & $u^{-1}tuv$ & $t[u,s]\otimes v + t\otimes tuvu^{-1}t^{-1}[s,u] + tuvu^{-1}[s,u]\otimes 1$ & Y\\\hline
$sv$ & $tu$ & $2tu[s,v]\otimes 1 + t\otimes tut^{-1}[s,v] + ut[v,s]\otimes 1 + t\otimes [v,s]u + 1\otimes [s,v]tu$ & Y\\\hline
$tu$ & $stuvu^{-1}t^{-1}$ & $0$ & Y\\\hline\hline
$a$ & $b$ & $-ab\otimes 1$ & N\\\hline
$a$ & $bc$ & $-acb\otimes 1$ & N\\\hline
$a$ & $bcd$ & $-acdb\otimes 1$ & N\\\hline
$ab$ & $ab$ & $-2ab\otimes ab + 2abab\otimes 1$ & Y\\\hline
$ab$ & $bc$ & $-bcba\otimes 1 + cbba\otimes 1 - ab\otimes bc$ & N\\\hline
$a$ & $aba^{-1}b^{-1}$ & $-2ba\otimes b^{-1} + b\otimes b^{-1}a + baa\otimes a^{-1}b^{-1}$ & Y\\\hline
\end{tblr}}
\caption{Some values of $\delta_2^{\sigma,\nabla\!_\mathcal{C}}$.}
\label{tab:2div}
\end{table}

\begin{example}\label{ex:second}
Some values of $\delta^{\sigma,\nabla\!_\mathcal{C}}_2$ are listed in Table \ref {tab:2div}. We put $a = \alpha_1$, $b = \beta_1$, $c = \alpha_2,\dotsc, s = \gamma_1$, $t = \gamma_2$\, and so on, and omit the vertical bar $|$ for the trace space for the sake of legibility. $[x,y] = xy-yx$ is the usual commutator.

It yields a non-zero value for some pairs $(x,y)$ of disjoint loops. Since the Goldman bracket and the Turaev cobracket vanish if loops have no (mutual or self-) intersections, so $\delta^{\sigma,\nabla\!_\mathcal{C}}_2$ cannot be written as a composite of them. There are hints of topological features, such as disjointness, reflected in the values (e.g., the sum of coefficients or the number of terms), but none are particularly clear. For similar arguments involving cyclic words, the reader may want to compare it with a combinatorial description of the Goldman--Turaev Lie bialgebra \cite{chas} and the Counting Intersection Theorem for the Thurston--Wolpert--Goldman Lie algebra \cite{CIT}.\\[-7pt]
\end{example}

The original motivation to define the maps $\Div_k^\nabla$ was to eliminate the choice of a connection in reformulating the divergence $\Div_1^\nabla$. Since we take the trace of $c_\nabla$, the off-diagonal entries of the connection form of $\nabla$ do not contribute to $\Div_1^\nabla$ at all, and the author expected that the conditions on $\Div_k^\nabla$, to which they certainly contribute, will restrict the choice of a connection. Therefore, we ask the following:

\begin{question}
Is there a (flat) connection $\nabla$ on $\Omega^1\mathbb{K}\pi$ that makes the operation $\delta_k^{\sigma,\nabla}$ (anti-)symmetric in its values?
\end{question}

\begin{remark} In the case of the tensor algebra $A = T(W)$, we have a canonical flat connection $\nabla\!_W$ explained below, and a similar statement holds, as we will see in the next section.\\
\end{remark}

\section{Props of Ribbon Graphs}\label{sec:ribbonop}

In this section, we will provide an explicit expression for $\Div^\nabla_k$ in the case of a tensor algebra and see a relation to operations associated with ribbon graphs introduced in \cite{ribbon} by Merkulov and Willwacher.

\begin{definition}\label{def:connW}
Let $W$ be a finite-dimensional $\mathbb{K}$-vector space and $T(W)$ the tensor algebra on $W$. The canonical flat  connection $\nabla\!_W$ on the free $T(W)^\mathrm{e}$-module $\Omega^1T(W)$ is defined by $\nabla\!_W(dw) = 0$ for all $w\in W$.
\end{definition}

\begin{definition}
Let $\pair{\cdot}{\cdot}\colon W\otimes W \to \mathbb{K}$ be a pairing. We extend the pairing to $T(W)\otimes T(W) \to T(W)\otimes T(W)$ by, for $u_i, w_j\in W$,
	\[
		\pair{u_1\cdots u_r}{w_1\cdots w_s} = \sum_{\substack{1\leq i\leq r\\1\leq j\leq s}} \pair {u_i}{w_j} w_1\cdots w_{j-1}u_{i+1}\cdots u_r\otimes u_1\cdots u_{i-1}w_{j+1}\cdots w_s\,.
	\]
	The associated \textit{Hamiltonian flow} is defined by 
	\[
		\mathrm{Ham}_{\pair{\cdot}{\cdot}}\colon |T(W)| \to \End_\mathbb{K}(T(W)): a \mapsto \pair{a}{\cdot}\,.
	\]
	We can show that the image is contained in $\Der_\mathbb{K}(T(W))$, and it is a Lie algebra morphism if the pairing is skew-symmetric.
\end{definition}

\begin{remark}
The map $\pair{\cdot}{\cdot}\colon T(W)\otimes T(W) \to T(W)\otimes T(W)$ satisfies the axiom of a \textit{double bracket} introduced in \cite{doublep} by Van den Bergh. In general, a double bracket $\Pi\colon A\otimes A\to A\otimes A$ induces the associated Hamiltonian flow $\mathrm{Ham}_\Pi\colon |A| \to \Der_\mathbb{K}(A)$. See \cite{toyo2} for the details.\\
\end{remark}

\begin{figure}[tb]
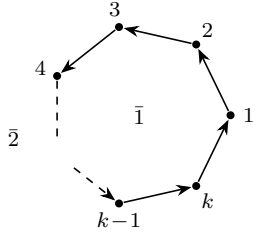

	\centering
	\RGloop
	\caption{The ribbon graph $L_k\in\mathcal{R}^k_{2,k}$. The labelling of the edges is suppressed, as any choice of these gives the same element in $\mathcal{RG}ra_1(2,k)$.}
	\label{fig:Lk}
\end{figure}

We have the following interpretation of the map $\delta^{\mathrm{Ham}_{\pair{\cdot}{\cdot}},\nabla\!_W}_k$ in terms of ribbon graphs:

\begin{theorem}\label{thm:ribbonop}
Suppose that the pairing $\pair{\cdot}{\cdot}\colon W\otimes W \to \mathbb{K}$ is skew-symmetric. Then, the map
\[
	(-1)^k\delta^{\mathrm{Ham}_{\pair{\cdot}{\cdot}},\nabla\!_W}_k\colon |T(W)|^{\otimes k} \to |T(W)|^{\otimes 2}
\]
coincides with the operation given by the ribbon graph $L_k$ in Figure \ref{fig:Lk} for each $k$.
\end{theorem}
\noindent The explanation and the proof of the proposition occupy the rest of the section.\\

Before the proof, we recall the definition of a ribbon graph and the construction of the associated operations. The notations below are aligned with the ones in \cite{ribbon}.

\begin{definition}\ 
\begin{itemize}
	\item A \textit{ribbon graph} is a tuple $\Gamma = (H,\tau_1,\tau_0,I)$ where $H$ is a non-empty finite set of \textit{half edges}, $\tau_1\colon H\to H$ is an involution (i.e., $\tau_1^2=\id$) with no fixed points, $\tau_0\colon H\to H$ is a bijection, and $I$ is a (possibly empty) finite set of \textit{isolated vertices}. The set of \textit{vertices}, \textit{edges} and \textit{boundary components} are defined by
	\[
		V(\Gamma) = (H/\tau_0)\sqcup I,\quad E(\Gamma) = H/\tau_1,\;\textrm{ and } \;B(\Gamma) = H/(\tau_1^{-1}\circ\tau_0)\sqcup I\,.
	\]
	\item For a vertex $v\in V(\Gamma)\!\setminus \!I$, let $H(v)\subset H$ be the orbit of $\tau_0$ representing $v$. This set is equipped with the natural cyclic order induced by $\tau_0$. The \textit{valency} of the vertex is the cardinality of the set $H(v)$. If $v\in I$, the valency is defined to be zero.
	\item A \textit{direction} of an edge $e = \{h,\tau_1(h)\}\in E(\Gamma)$ is a total order on the set $\{h,\tau_1(h)\}$. A \textit{directed} ribbon graph is a ribbon graph with a choice of direction for each edge.
	\item A \textit{labelling} of a ribbon graph $\Gamma$ with $n$ vertices, $\ell$ edges and $m$ boundary components is a tuple of bijections
\begin{align*}
	V(\Gamma) \to \{1,\dotsc,n\},\; E(\Gamma) \to \{1,\dotsc,\ell\} \textrm{ and } B(\Gamma) \to \{\bar 1,\dotsc, \bar m\}\,.
\end{align*}
	\item $\mathcal R_{m,n}^\ell$ is the set of (isomorphism classes of) all such directed labelled ribbon graphs.
\end{itemize}
\end{definition}

A ribbon graph is depicted as a planar graph as in Figure \ref{fig:Lk}, with the cyclic order at each vertex specified in counterclockwise rotation.\\

Finally, the prop $\mathcal{RG}ra_1$ is defined as follows: as a $\mathbb{K}$-vector space,
\begin{align*}
	\mathcal{RG}ra_1(m,n) = \bigoplus_{\ell\geq 0} \mathbb{K}\mathcal R_{m,n}^\ell/(S_\ell\rtimes S_2^{\times \ell})\\[-20pt]
\end{align*}
for $m,n\geq 1$, where the symmetric group $S_\ell$ acts on it by the re-numbering of edges, and each copy of $S_2$ acts by the corresponding edge flip with the factor $(-1)$ for each transposition. The prop structure is given by a composition of ribbon graphs, which is an analogue of surface glueing; we will not use the explicit formula below and refer to Section 4.2 of \cite{ribbon} for the detail.\\

The prop $\mathcal{RG}ra_1$ has a canonical representation on the space of cyclic words $|T(W)|$ over a vector space $W$ with a skew-symmetric pairing (Theorem 4.2.2 of \cite{ribbon}), and we shall briefly recall it. With a ribbon graph $\Gamma\in\mathcal{R}^\ell_{m,n}$ is associated a $\mathbb{K}$-linear map $\Gamma\colon |T(W)|^{\otimes n}\to |T(W)|^{\otimes m}$ by the following procedure. To compute $\Gamma(\mathbf{w}_1,\dotsc,\mathbf{w}_n)$ for cyclic words $\mathbf{w}_i = w^{(i)}_1\cdots w^{(i)}_{r_i}$ with $w^{(i)}_j\in W$:
\begin{enumerate}[(1)]
	\item Take a cyclic-order preserving injection $q_i\colon H(i)\to \{1,\dotsc,r_i\}$ for each $i\in V(\Gamma)$. Here each vertex $v$ is identified with its label $i$, and the set $\{1,\dotsc,r_i\}$ is equipped with the cyclic order $1\leq 2\leq \cdots\leq r_i\leq 1$. This gives a cyclic arrangement of $w^{(i)}_j$'s around the vertex $i$ with $w^{(i)}_{q_i(h)}$ \textit{on} the half-edge $h$ and other letters \textit{in-between}. If not possible (i.e., if $ \# H(i)>r_i$), the value of $\Gamma(\mathbf{w}_1,\dotsc,\mathbf{w}_n)$ is zero.
	\item For each directed edge $e = \{h,\tau_1(h):h\leq \tau_1(h)\}$ from a vertex $i$ to the other $i'$ in this direction, take a pairing $\pair{w^{(i)}_{q_i(h)}}{w^{(i')}_{q_{i'}(\tau_1(h))}}$. Denote by $\lambda$ the product of these scalars over all edges. The letters used in the pairing are thus \textit{consumed}.
	\item Obtain the cyclic words $\mathbf{w}'_{k}$ by collecting the surviving letters $w^{(i)}_j$ along the $k$-th boundary of $\Gamma$.
	\item $\Gamma(\mathbf{w}_1,\dotsc,\mathbf{w}_n)$ is the sum of $\lambda \cdot\mathbf{w}'_1\otimes \cdots\otimes \mathbf{w}'_m$ for all possible combination of injections $q = (q_i)_{1\leq i\leq n}$.
\end{enumerate}
Then, this only depends on the class in $\mathcal{RG}ra_1(m,n)$ by the skew-symmetry of the pairing.

Moreover, there is also a canonical map $\mathcal{L}ieb^\diamond\to\mathcal{RG}ra_1$ from the prop governing involutive Lie bialgebras given by
\begin{align*}
	(\textrm{Lie bracket}) \mapsto \RGbar\;\textrm{ and }\; (\textrm{Lie cobracket}) \mapsto \RGtadpole = L_1\,.
\end{align*}
For the proof, see Theorem 4.2.3 of \cite{ribbon}. Therefore, an $\mathcal{RG}ra_1$-algebra can be considered as an enhancement of an involutive Lie bialgebra, and $|T(W)|$ is automatically endowed with an involutive Lie bialgebra structure. In general, for each integer $d$, we can define a $\mathbb{Z}$-graded prop $\mathcal{RG}ra_d$ with the canonical representation on a graded vector space with graded-symmetric pairing.\\

\noindent\textbf{Proof of Theorem \ref{thm:ribbonop}.} The proof is done by direct computation. Take $\mathbf{w}_1,\dotsc,\mathbf{w}_k\in |T(W)|$ of the form $\mathbf{w}_i = w^{(i)}_1\cdots w^{(i)}_{r_i}$ with $w^{(i)}_j\in W$. For $w\in W$, we have 
\begin{align*}
	\mathrm{Ham}_{\pair{\cdot}{\cdot}}(\mathbf{w}_i)(w) &= \sum_{1\leq s_i\leq r_i} \pair{w^{(i)}_{s_i}}{w}w^{(i)}_{s_i+1}\cdots w^{(i)}_{s_i-1} ,\\
	c_{\nabla\!_W}(\mathrm{Ham}_{\pair{\cdot}{\cdot}}(\mathbf{w}_i))(dw) &= -d\Bigg(\sum_{1\leq s_i\leq r_i} \pair{w^{(i)}_{s_i}}{w}w^{(i)}_{s_i+1}\cdots w^{(i)}_{s_i-1} \Bigg)\\
	&= -\sum_{\substack{1\leq s_i\leq r_i\\1\leq t_i\leq r_i\\s_i\neq t_i}} \pair{w^{(i)}_{s_i}}{w}w^{(i)}_{s_i+1}\cdots w^{(i)}_{t_i-1} dw^{(i)}_{t_i} w^{(i)}_{t_i+1}\cdots w^{(i)}_{s_i-1} ,
\end{align*}
\begin{align*}
	&c_{\nabla\!_W}(\mathrm{Ham}_{\pair{\cdot}{\cdot}}(\mathbf{w}_1))\circ\cdots\circ c_{\nabla\!_W}(\mathrm{Ham}_{\pair{\cdot}{\cdot}}(\mathbf{w}_k))(dw)\\
	&= c_{\nabla\!_W}(\mathrm{Ham}_{\pair{\cdot}{\cdot}}(\mathbf{w}_1))\circ\cdots\circ c_{\nabla\!_W}(\mathrm{Ham}_{\pair{\cdot}{\cdot}}(\mathbf{w}_{k-1}))\Bigg( -\sum_{\substack{1\leq s_i\leq r_i\\1\leq t_i\leq r_i\\s_i\neq t_i}} \pair{w^{(k)}_{s_k}}{w}w^{(k)}_{s_k+1}\cdots w^{(k)}_{t_k-1} dw^{(k)}_{t_k} w^{(k)}_{t_k+1}\cdots w^{(k)}_{s_k-1} \Bigg)\\
	&= \cdots\textrm{(repeatedly applied)}\cdots\\[5pt]
	&= \sum_{1\leq i\leq k}\sum_{\substack{1\leq s_i\leq r_i\\1\leq t_i\leq r_i\\s_i\neq t_i}} \pair{w^{(1)}_{s_1}}{w^{(2)}_{t_2}}\cdots \pair{w^{(k)}_{s_k}}{w} (-1)^k\\
	&\qquad w^{(k)}_{s_k+1}\cdots w^{(k)}_{t_k-1}\cdots w^{(2)}_{s_2+1}\cdots w^{(2)}_{t_2-1}w^{(1)}_{s_1+1}\cdots w^{(1)}_{t_1-1} dw^{(1)}_{t_1} w^{(1)}_{t_1+1}\cdots w^{(1)}_{s_1-1}w^{(2)}_{t_2+1}\cdots w^{(2)}_{s_2-1}\cdots w^{(k)}_{t_k+1}\cdots w^{(k)}_{s_k-1},
\end{align*}
and, by taking the trace, we have
\begin{align*}
	(-1)^k\delta^{\mathrm{Ham}_{\pair{\cdot}{\cdot}},\nabla\!_W}_k(\mathbf{w}_1,\dotsc,\mathbf{w}_k) = \sum_{1\leq i\leq k}\sum_{\substack{1\leq s_i\leq r_i\\1\leq t_i\leq r_i\\s_i\neq t_i}}& \pair{w^{(1)}_{s_1}}{w^{(2)}_{t_2}}\cdots \pair{w^{(k)}_{s_k}}{w^{(1)}_{t_1}}\\[-20pt]
	&\quad|w^{(k)}_{s_k+1}\cdots w^{(k)}_{t_k-1}\cdots w^{(2)}_{s_2+1}\cdots w^{(2)}_{t_2-1}w^{(1)}_{s_1+1}\cdots w^{(1)}_{t_1-1}|\\
	&\hspace{20pt}\otimes |w^{(1)}_{t_1+1}\cdots w^{(1)}_{s_1-1}w^{(2)}_{t_2+1}\cdots w^{(2)}_{s_2-1}\cdots w^{(k)}_{t_k+1}\cdots w^{(k)}_{s_k-1}|\,.
\end{align*}
which is exactly equal to $L_k(\mathbf{w}_1,\dotsc,\mathbf{w}_k)$ by the definition of the ribbon graph operation. \qed\\

\begin{remark}\label{rem:unsuc}
The original motivation (in addition to giving a constraint on the choice of connections) was to craft a family of loop operations, each corresponding to an element of $\mathcal{RG}ra_1$ being an extension of the Goldman--Turaev Lie bialgebra. The operations given by the graphs
\[
	\RGbar \textrm{ and }\RGtadpole
\]
are the associated graded of the Goldman bracket and the Turaev cobracket, respectively, for the surface with only one boundary component. See the survey \cite{gtjohnson} by Kawazumi--Kuno for the details.

There had already been some attempts by several people (and the author is, again, thankful to Sergei Merkulov for sharing their past attempts on it), but all were unfortunately unsuccessful; operations given by the smoothing-intersections type argument are not well-defined when two or more intersections are resolved at the same time. In contrast, our definition of the operations $\delta^{\psi, \nabla}_k$ are purely algebraic with a lack of topological perspective. The interpretation of these operations in terms of the topology of curves would be interesting. 

Another result in the paper \cite{ribbon} is that there is a map of props $\mathcal Ho\mathcal Lieb_{\hspace{0.5pt}0,1} \to \mathcal{RG}ra_1$ (Theorem 4.3.3 of \cite{ribbon}). It would also be interesting to investigate the $\mathcal{H}o\mathcal{L}ieb_{0,1}$-relations in the case of the group algebra $\mathbb{K}\pi$.\\
\end{remark}

\section{Non-commutative de Rham Cohomology and Vector Bundles}\label{sec:tdr}

The quantity $c_\nabla$ defined in Section \ref{sec:higher} was essential for the construction. To see what it really is, we will develop a theory on a non-commutative version of de Rham cohomology with a vector bundle coefficient. We will relate the space of endomorphism-valued non-commutative differential forms to the Chevalley--Eilenberg (CE) complex over a derivation Lie algebra in Theorem \ref{thm:twdr} and show that the alternated version $\Div^\nabla_k\circ\,\mathrm{alt}$ is a $k$-cocycle in the CE complex.\\[-7pt]

Let $B$ be a unital associative $\mathbb{K}$-algebra.

\begin{definition}
The \textit{graded} trace space $\DR^\bullet\!B = \Omega^\bullet B/[\Omega^\bullet B,\Omega^\bullet B]$ of $\Omega^\bullet B$ is the space of \textit{de Rham differential forms}. Being a graded derivation, the exterior derivative $d$ descends to $\DR^\bullet\!B$.
\end{definition}

Note that an element of $\DR^n\!B$ can be written in the form $|b_0db_1\cdots db_n|$ with $b_i\in B$.
\begin{definition}
For a $B$-module $M$, we define the space of \textit{endomorphism-valued differential forms} by
\[
	\Omega^\bullet(B,\End M) = \Hom_B(M,\Omega^\bullet B\otimes_BM).
\]
For a connection $\nabla$ on $M$, we define the covariant derivative $D_\nabla$ on the space $\Omega^\bullet(B,\End M)$ by
\begin{align*}
	D_\nabla:\Omega^\bullet(B,\End M)&\to \Omega^{\bullet+1}(B,\End M)\\
	\xi &\mapsto [\nabla,\xi] := \nabla\circ\xi - (-1)^{|\xi|} (\id\otimes \,\xi)\circ\nabla .
\end{align*}
Here we used the extension of $\nabla$ to $\Omega^\bullet B\otimes_BM$ defined in Section \ref{sec:higher}.
\end{definition}

\begin{remark}
The notation $\Omega^\bullet(B,\End M)$ is purely formal. This is not equal to the $B$-tensor product of $\Omega^\bullet B$ and $\End_B(M)$: The latter does not even have a $B$-module structure unless $B$ is commutative. When $M$ is dualisable, we have a natural isomorphism $\Omega^\bullet(B,\End M)\cong M^*\otimes_B \Omega^\bullet B\otimes_BM$.

Also, since we use the Koszul sign convention, we understand $(\id\otimes\,\xi)\circ\eta$ as $(-1)^{|\xi||\eta|}\eta'\otimes (\xi\circ\eta'')$ for $\xi,\eta \in \Omega^\bullet(B,\End M)$. It is similar for $\nabla$, which is not $B$-linear.
\end{remark}

Now recall that a \textit{curved differential graded Lie algebra} is a triple $(\mathscr{L}, d_\mathscr{L}, R_\mathscr{L})$ where $\mathscr{L}$ is a graded Lie algebra (with Koszul signs), $d_\mathscr{L}$ is a degree $+1$ Lie algebra derivation, and $R_\mathscr{L}$ is a degree $+2$ element in $\mathscr{L}$ such that $d_\mathscr{L}^2 = [R_\mathscr{L},\cdot\,]$.

\begin{proposition}
$(\Omega^\bullet(B,\End M), D_\nabla, \nabla^2)$ is a curved differential graded Lie algebra with the graded commutator $[\xi,\eta] = (\id\otimes\,\xi)\circ\eta -(-1)^{|\xi||\eta|}(\id\otimes\,\eta)\circ\xi$.
\end{proposition}
\noindent Proof. Since the Lie bracket is induced from the graded associative algebra structure $\xi\cdot \eta = (\id\otimes\,\xi)\circ\eta$ on $\Omega^\bullet(B,\End M)$, it is automatically a graded Lie algebra. Next, the covariant derivative is clearly a derivation of degree $+1$. Finally, since the curvature $\nabla^2\colon M \to \Omega^2B\otimes_BM$ is $B$-linear, we have $\nabla^2 \in \Omega^2(B,\End M)$ and $D_\nabla^2 = [\nabla^2,\cdot\,]$ by definition. \qed

\begin{definition}
Let $M$ be $B$-dualisable. The \textit{trace map} on $\Omega^\bullet(B,\End M)$ with values in the non-commutative de Rham complex $\DR^\bullet\!B$ is defined by
\begin{align*}
	\Tr:\Omega^\bullet(B,\End M)&\to \DR^\bullet\!B\\
	\theta\otimes\omega\otimes m &\mapsto |\theta(m)\omega|
\end{align*}
using the natural surjection $|\Omega^\bullet B| \twoheadrightarrow \DR^\bullet\! B$.
\end{definition}

\begin{proposition}\label{prop:HTtraceDR}
The trace map
\[
	\Tr:(\Omega^\bullet(B,\End M),D_\nabla)\to (\DR^\bullet\!B,d)
\]
is a curved differential graded Lie algebra homomorphism. Here $\DR^\bullet \!B$ is regarded as an abelian Lie algebra.
\end{proposition}
\noindent Proof. Putting $\nabla = \nabla'\otimes\nabla'':M\to \Omega^1B\otimes_BM$, one has
\begin{align*}
	D_\nabla(\theta\otimes\omega\otimes m) &= [\nabla,\theta\otimes\omega\otimes m]\\
	&= d\theta(\cdot)\omega\otimes m + \theta(\cdot)d\omega\otimes m + (-1)^{|\omega|}\theta(\cdot)\omega\nabla'(m)\otimes\nabla''(m)\\
	& \qquad  - (-1)^{|\omega|+|\omega|} \nabla'(\cdot)\theta(\nabla''(\cdot))\omega\otimes m\,,\\
	\intertext{so that}
	\Tr D_\nabla(\theta\otimes\omega\otimes m) &= |d\theta(m)\omega+ \theta(m)d\omega + (-1)^{|\omega|}\theta(\nabla''(m))\omega\nabla'(m) -  \nabla'(m)\theta(\nabla''(m))\omega|\\
	&=  d|\theta(m)\omega| = d\Tr(\theta\otimes\omega\otimes m),
\end{align*}
which shows that $\Tr$ intertwines the differentials. The proof that the trace map is a Lie algebra homomorphism is analogous to the proof that the trace map is cyclic in its argument. \qed\\

Since de Rham cohomology is closely related to Chevalley--Eilenberg cohomology in the usual differential geometry, we shall now relate these in a non-commutative setting.

\begin{definition}\ 
\begin{itemize}
	\item For a Lie algebra $\mathfrak{g}$ over $\mathbb{K}$ and a $\mathfrak{g}$-module $V$, $C^n_\mathrm{CE}(\mathfrak{g},V) = \Hom_\mathbb{K}(\wedge^n \mathfrak{g},V)$ is the Chevalley-Eilenberg complex with the differential
	\begin{align*}
		(d_\mathrm{CE}\psi)(x_0,\dotsc,x_n) &= \sum_{0\leq i\leq n} (-1)^i x_i\cdot \psi(x_0,\dotsc,\check x_i,\dotsc,x_n)\\ 
		&\qquad\qquad  + \sum_{0\leq i<j\leq n} (-1)^{i+j} \psi([x_i,x_j],x_0,\dotsc,\check x_i,\dotsc,\check x_j,\dotsc,x_n).
	\end{align*}
	If $V$ is an associative $\mathbb{K}$-algebra with a $\mathfrak{g}$-equivariant multiplication, the complex is a differential graded algebra with the shuffle product.
	\item For a derivation action $(\mathfrak d, \varphi,\rho)$ on a $B$-module $M$, the space $\End_B(M)$ is regarded as a $\mathfrak d$-module by $f\cdot\mu = [\rho(f),\mu]$ for $f\in\mathfrak d$ and $\mu\in\End_B(M)$. Here $[\rho(f),\mu]$ is again $B$-linear by the compatibility condition for a derivation action.
	\item Given a derivation action $(\mathfrak d, \varphi,\rho)$ on $M$, the maps
	\[
		\iota\colon \Omega^\bullet(B,\End M) \to C^\bullet_\mathrm{CE}(\mathfrak d,\End_B(M))\quad\textrm{and}\quad\iota\colon \DR^\bullet\!B \to C^\bullet_\mathrm{CE}(\mathfrak d,|B|)
	\]
	are defined by, for $f_i\in\mathfrak d$ and $\alpha\in \Omega^n(B,\End M)$ or $\DR^n\!B$,
	\[
		\iota(\alpha)(f_1\wedge\cdots\wedge f_n) = i_{f_n}\!\cdots i_{f_1}\alpha\,.
	\]
	Here $i_{\varphi(f)}\otimes \id$ is briefly denoted by $i_f$.
\end{itemize}
\end{definition}

\begin{remark}
The map $c_\nabla$ defined in Section \ref{sec:higher} lives in $C^1_\mathrm{CE}(\mathfrak d,\End_B(M))$.
\end{remark}

The following is the key property of the map $c_\nabla$.

\begin{proposition}\label{prop:connMC}
If $\nabla$ is flat, $c_\nabla$ is a Maurer--Cartan element with respect to the differential $d_\mathrm{CE}$. Therefore, $d_\mathrm{CE}+[c_\nabla,\cdot\,]$ is also a differential.
\end{proposition}
\noindent Proof. Let $R = \nabla^2$ be the curvature. For a general connection $\nabla$, we have
\begin{align*}
	(d_\mathrm{CE}c_\nabla)(f\wedge g) + [c_\nabla(f),c_\nabla(g)] = \iota(R)(f\wedge g).
\end{align*}
The proof is essentially the same as Proposition 4.11 in \cite{toyo}. By the definition of the shuffle product, $[c_\nabla(f),c_\nabla(g)]$ is equal to $(c_\nabla\circ c_\nabla)(f\wedge g)$; this implies $d_\mathrm{CE}c_\nabla + c_\nabla\circ c_\nabla = \iota(R)$, which is zero if $\nabla$ is flat. \qed

\begin{theorem}\label{thm:twdr}
The following is a commutative diagram of curved complexes:
\[
\begin{tikzcd}[cramped]
	\big(\Omega^\bullet(B,\End M),D_\nabla\big) \arrow[r, "\Tr"] \arrow[d, "\iota"] & (\DR^\bullet\!B,d) \arrow[d, "\iota"]\\
	\big(C^\bullet_\mathrm{CE}(\mathfrak d,\End_B(M)),\,d_\mathrm{CE}+[c_\nabla,\cdot\,]\big) \arrow[r, "\Tr"] & \big(C^\bullet_\mathrm{CE}(\mathfrak d,|B|),d_\mathrm{CE}\big)\,.\\[-10pt]
\end{tikzcd}
\]
\end{theorem}

We need some preparations for the proof.
\begin{proposition}\label{prop:iota}
Both $\iota$'s are cochain maps.
\end{proposition}
\noindent Proof. We only check the left one since the right one is easier. For $\theta\otimes\omega\otimes m\in M^*\otimes_B\Omega^nB\otimes_BM$, we have (compare with the computation in the proof of Proposition \ref{prop:HTtraceDR}), denoting $i_{\varphi(f)}$ simply by $i_f$,
\begin{align*}
	\iota(&D_\nabla(\theta\otimes\omega\otimes m))(f_0\wedge\cdots\wedge f_n)\\
	&= i_{f_n}\cdots i_{f_0}(d\theta\,\omega\otimes m + \theta d\omega\otimes m + (-1)^n\theta\omega\nabla m - (\theta\circ\nabla)\cdot\omega\otimes m)\\
	&= \sum_{0\leq j\leq n}(-1)^j (f_j\circ\theta)\cdot i_{f_n}\cdots\check{i_{f_j}}\cdots i_{f_0}\omega\cdot m \:+\: \theta \cdot i_{f_n}\cdots i_{f_0}d\omega\cdot m\\
	&\qquad+ \sum_{0\leq j\leq n} \left((-1)^{n+(n-j)}\theta\cdot i_{f_n}\cdots\check{i_{f_j}}\cdots i_{f_0}\omega\cdot i_{f_j}\nabla m - (-1)^j (\theta\circ i_{f_j}\nabla)\cdot i_{f_n}\cdots\check{i_{f_j}}\cdots i_{f_0}\omega\otimes m\right).
\end{align*}
On the other hand, we have
\begin{align*}
	((d&_\mathrm{CE} + [c_\nabla,\cdot\,])\circ\,\iota)(\theta\otimes\omega\otimes m))(f_0\wedge\cdots\wedge f_n)\\
	&= \sum_{0\leq j\leq n}(-1)^j [f_j,\theta\cdot i_{f_n}\cdots\check{i_{f_j}}\cdots i_{f_0}\omega\cdot m] + \sum_{0\leq j<k\leq n} (-1)^{j+k} \theta\cdot i_{f_n}\cdots\check{i_{f_k}}\cdots\check{i_{f_j}}\cdots i_{f_0}i_{[f_j,f_k]}\omega\cdot m\\
	&\qquad + \sum_{0\leq j\leq n} (-1)^j [c_\nabla(f_j),\theta\cdot i_{f_n}\cdots\check{i_{f_j}}\cdots i_{f_0}\omega\cdot m]\\
	&=  \sum_{0\leq j\leq n}(-1)^j [i_{f_j}\nabla,\theta\cdot i_{f_n}\cdots\check{i_{f_j}}\cdots i_{f_0}\omega\cdot m] + \sum_{0\leq j<k\leq n} (-1)^{j+k} \theta\cdot i_{f_n}\cdots\check{i_{f_k}}\cdots\check{i_{f_j}}\cdots i_{f_0}i_{[f_j,f_k]}\omega\cdot m\\
	&= \sum_{0\leq j\leq n}(-1)^j \Big((f_j\circ\theta)\cdot i_{f_n}\cdots\check{i_{f_j}}\cdots i_{f_0}\omega\cdot m + \theta\cdot L_{f_j}(i_{f_n}\cdots\check{i_{f_j}}\cdots i_{f_0}\omega)\cdot m \\
	&\hspace{100pt} + \theta\cdot i_{f_n}\cdots\check{i_{f_j}}\cdots i_{f_0}\omega\cdot i_{f_j}\nabla m - (\theta\circ i_{f_j}\nabla)\cdot i_{f_n}\cdots\check{i_{f_j}}\cdots i_{f_0}\omega\cdot m\Big)\\
	&\qquad + \sum_{0\leq j<k\leq n} (-1)^{j+k} \theta\cdot i_{f_n}\cdots\check{i_{f_k}}\cdots\check{i_{f_j}}\cdots i_{f_0}i_{[f_j,f_k]}\omega\cdot m.
\end{align*}
Finally, we are done if we prove, on $\Omega^nB$,
\[
	i_{f_n}\cdots i_{f_0}d = \sum_{0\leq j\leq n}(-1)^j L_{f_j}i_{f_n}\cdots\check{i_{f_j}}\cdots i_{f_0} +  \sum_{0\leq j<k\leq n} (-1)^{j+k} i_{f_n}\cdots\check{i_{f_k}}\cdots\check{i_{f_j}}\cdots i_{f_0}i_{[f_j,f_k]}.
\]
We proceed by induction. If $n=0$, this is equivalent to $i_{f_0}d = L_{f_0}$, which is true because $\Omega^0B$ is simply equal to $B$. For $n\geq 1$, by the induction hypothesis,
\begin{align*}
	i_{f_n}\cdots i_{f_0}d &=i_{f_n}\cdots i_{f_1}L_{f_0} - (i_{f_n}\cdots i_{f_1}d)i_{f_0}\\
	&= L_{f_0}i_{f_n}\cdots i_{f_1} + \sum_{1\leq k\leq n} i_{f_n}\cdots i_{[f_k,f_0]}\cdots i_{f_1}\\
	&\qquad - \left( \sum_{1\leq j\leq n}(-1)^{j-1} L_{f_j}i_{f_n}\cdots\check{i_{f_j}}\cdots i_{f_1} +  \sum_{1\leq j<k\leq n} (-1)^{j+k} i_{f_n}\cdots\check{i_{f_k}}\cdots\check{i_{f_j}}\cdots i_{f_1}i_{[f_j,f_k]}  \right)i_{f_0}\\
	&=  \sum_{0\leq j\leq n}(-1)^j L_{f_j}i_{f_n}\cdots\check{i_{f_j}}\cdots i_{f_0} +  \sum_{0\leq j<k\leq n} (-1)^{j+k} i_{f_n}\cdots\check{i_{f_k}}\cdots\check{i_{f_j}}\cdots i_{f_0}i_{[f_j,f_k]}.
\end{align*}
This completes the proof.\qed\\

\noindent\textbf{Proof of Theorem \ref{thm:twdr}.} The square is commutative since contractions and trace maps are mutually compatible. Furthermore, $\Tr:\End_B(M)\to|B|$ is automatically a $\Der_\mathbb{K}(B)$-module map; hence is a $\mathfrak d$-module map by restriction. Therefore, $\Tr$ commutes with $d_\mathrm{CE}$. For any $\psi\in C^\bullet_\mathrm{CE}(\mathfrak d,\End_B(M))$, we have
\[
	\Tr(d_\mathrm{CE}(\psi) + [c_\nabla,\psi]) = \Tr(d_\mathrm{CE}(\psi)) = d_\mathrm{CE}\Tr(\psi).
\]
by the cyclic symmetry of the trace map. This shows that $\Tr$ in the bottom row is a cochain map. The top row is a cochain map by Proposition \ref{prop:HTtraceDR}. Finally, Proposition \ref{prop:iota} ensures that $\iota$'s are cochain maps.\qed\\

We will now look at the map $\Div^\nabla_k$ with these tools available. To apply the above, we only consider the skew-symmetrised version $\Div^\nabla_k\circ\,\mathrm{alt}$, where $\mathrm{alt}$ is defined by 
\[
	\mathrm{alt}(f_1\wedge\dotsc\wedge f_k) = \sum_{s\in S_k} \mathrm{sgn}(s) f_{s(1)}\otimes \cdots\otimes f_{s(k)}\,.
\]

\begin{lemma}
The map $\Div^\nabla_k\!\circ\,\mathrm{alt}\colon \wedge^k \!\mathfrak d \to |B|$ is equal to $\Tr(c_\nabla{}^k)\in C^k_\mathrm{CE}(\mathfrak d, |B|)$. Here, $c_\nabla{}^k$ is the $k$-fold product in the CE complex.
\end{lemma}
\noindent Proof. This follows from the definition of $\Div^\nabla_k$ and the shuffle product. \qed\\[-5pt]

\begin{proposition}\label{prop:cocyczero}
For $k$ even, $\Div^\nabla_k\!\circ\,\mathrm{alt} = \Tr(c_\nabla{}^k)$ is zero.
\end{proposition}
\noindent Proof. By the cyclic symmetry of the trace, we have
\[
	\Tr(c_\nabla{}^{k}) = \Tr(c_\nabla \circ c_\nabla{}^{k-1}) = (-1)^{1\cdot (k-1)} \Tr(c_\nabla{}^{k-1} \circ c_\nabla) = -\Tr(c_\nabla{}^{k})\,.
\]
by the Koszul sign convention. This completes the proof. \qed\\[-5pt]

\begin{theorem}\label{thm:kcocyc}
If $\nabla$ is flat, the map $\Div^\nabla_k\!\circ\,\mathrm{alt}\in C^k_\mathrm{CE}(\mathfrak d, |B|)$ is a Lie algebra $k$-cocycle. If $k$ is even, the cocycle vanishes.
\end{theorem}
\noindent Proof. If $k$ is even, the map $\Div^\nabla_k\!\circ\,\mathrm{alt}$ itself is zero by Proposition \ref{prop:cocyczero}. Suppose $k$ is odd. Since $\nabla$ is flat, we have $(d_\mathrm{CE}+[c_\nabla,\cdot\,])(c_\nabla) = c_\nabla{}^2$ by Proposition \ref{prop:connMC}. In addition, $d_\mathrm{CE}+[c_\nabla,\cdot\,]$ is a degree $+1$ derivation of the CE complex, which is a differential graded algebra. Therefore, we have
\begin{align*}
	(d_\mathrm{CE}+[c_\nabla,\cdot\,])(c_\nabla{}^k) = \sum_{1\leq i\leq k}(-1)^{i-1} c_\nabla{}^{(i-1) + 2 + (k-i)} = c_\nabla{}^{k+1}.
\end{align*}
Taking the trace, we obtain $d_\mathrm{CE}(\Tr(c_\nabla{}^k)) = \Tr(c_\nabla{}^{k+1})$ in $C^{k+1}_\mathrm{CE}(\mathfrak d, |B|)$ by Theorem \ref{thm:twdr}. This is equal to zero by the corollary above.\qed\\

In the case of the tensor algebra $A = T(W)$ over a finite-dimensional $\mathbb{K}$-vector space $W$, we have the following non-vanishing theorem. Recall the result of Fuks (Theorem 2.1.6 of \cite{fuks}) which states that the cohomology ring $H_\mathrm{CE}^\bullet(\mathfrak{gl}(W),\mathbb{K})$ is isomorphic to the exterior algebra generated by $\varphi_k$ for $1\leq k\leq 2\dim(W)-1$ and $k$ odd, where
\begin{align*}
	\varphi_k\colon \wedge^k \mathfrak{gl}(W) &\to \mathbb{K}\\
	A_1\wedge\cdots\wedge A_k &\mapsto \sum_{s\in S_k}\mathrm{sgn}(s) \Tr(A_{s(1)}\cdots A_{s(k)}).
\end{align*}

\begin{theorem}\label{thm:dercohom}
The cohomology class of $\Tr(c_{\nabla\!_W}^{\;\;\;\,k})$ does not vanish in $H_\mathrm{CE}^k(\Der_\mathbb{K}(T(W)),|T(W)^\mathrm{e}|)$ for odd $k$.
\end{theorem}
\noindent Proof. Put $D = \Der_\mathbb{K}(T(W))$ and $V = |T(W)^\mathrm{e}|$. These spaces are naturally graded, and the $D$-action on $V$ has degree 0. Let $D^{(0)}$ be the Lie subalgebra of derivations with degree 0, and $j\colon D^{(0)} \hookrightarrow D$ be the natural inclusion, which is a Lie algebra homomorphism. We have an identification $D^{(0)} = \mathfrak{gl}(W)$. The pull-back map reads
\[
	H_\mathrm{CE}^\bullet(D,V) \xrightarrow{j^*} H_\mathrm{CE}^\bullet(\mathfrak{gl}(W),j^*V),
\]
where $j^*V$ is the restriction of $V$ by $j$. 

Now consider the Euler operator $\mathsf{eu}$, which is a unique degree $0$ derivation on $T(W)$ specified by $\mathsf{eu}(w) = w$ for all $w\in W$. It generates the centre of $\mathfrak{gl}(W)$ and acts on $V$ by multiplication-by-degree. Then, the ``centre kills'' argument shows that $H_\mathrm{CE}^\bullet(\mathfrak{gl}(W),j^*V) \cong H_\mathrm{CE}^\bullet(\mathfrak{gl}(W),\mathbb{K})$ since $j^*V^{(0)} = \mathbb{K}$. 

On the other hand, for $f\in D^{(0)}$, we have
\[
	(c_{\nabla\!_W}(f))(dw) = i_f\nabla\!_W(dw) - L_f (dw) = -df(w).
\]
Since $f$ is just a matrix over $\mathbb{K}$, the endomorphism $c_{\nabla\!_W}(f)\in \End_{T(W)^\mathrm{e}}\Omega^1T(W)$ is represented by the matrix $-f$. Therefore, we have
\[
	\Tr(c_{\nabla\!_W}^{\;\;\;k})(f_1\wedge\cdots\wedge f_k) = (-1)^k\varphi_k(f_1\wedge\cdots\wedge f_k)
\]
for $f_i\in D^{(0)}$, which implies $j^*(\Tr(c_{\nabla\!_W}^{\;\;\;k})) =  (-1)^k\varphi_k$. By the above result of Fuks, this does not vanish, and so does $\Tr(c_{\nabla\!_W}^{\;\;\;k})$. \qed\\

For many classical results on the commutative counterpart, see Chapter 2 of the book \cite{fuks} by Fuks. For the cohomology of the derivation Lie algebra on a free algebra over an (augmented) operad with various other coefficients, see the paper \cite{dotsenko} by Dotsenko.

\begin{remark}
In the case of the tensor algebra $T(W)$, an alternative proof of the map $\delta^{\mathrm{Ham}_{\pair{\cdot}{\cdot}},\nabla\!_W}_k$ in Proposition \ref{thm:ribbonop} being a Lie algebra cocycle for $k\equiv 1\!\!\mod 4$ is available in the framework of the ribbon graph complex by Merkulov and Willwacher \cite{ribbon}. First, using the notation there, the ribbon graph complex is given as the deformation complex
\[
	\mathsf{RGC}_{1,1} := \mathsf{Def}(\mathcal Lieb \xrightarrow{s_1} \mathcal{RG}ra_1)
\]
associated with the properad morphism $s_1$ given by
\begin{align*}
	(\textrm{Lie bracket}) \mapsto \RGbar\;\textrm{ and }\; (\textrm{Lie cobracket}) \mapsto 0\,.
\end{align*}
This complex is automatically a differential graded Lie algebra with differential induced by $s_1$, and it is naturally isomorphic to the space of all infinite series of directed labelled ribbon graphs with their vertex- and boundary-labelling skew-symmetrised. Furthermore, the differential is pictorially given by the vertex expansion, which increases the number of edges by exactly one and preserves bivalent vertices. In addition, the element corresponding to the graph $L_k$ is non-zero if and only if $k\equiv 1\!\!\mod 4$. From these facts, it is readily seen that the differential of $L_k$ is always zero.

Finally, the functorial property of the deformation complex gives the push-out map
\[
	\mathsf{RGC}_{1,1} = \mathsf{Def}(\mathcal Lieb \xrightarrow{s_1} \mathcal{RG}ra_1) \to \mathsf{Def}(\mathcal Lieb \to \mathcal End_{|T(W)|}) = \hat C^\bullet_\mathrm{CE}(|T(W)|, \wedge |T(W)|)\,,
\]
which is a differential graded Lie algebra morphism. Here, the CE complex in the last term is suitably completed. Under this morphism, the element $L_k$ is sent to $\mathrm{alt}\circ \delta^{\mathrm{Ham}_{\pair{\cdot}{\cdot}},\nabla\!_W}_k$, which ensures its cocycle property. When $k\equiv 1\!\!\mod 4$, $L_k$ is already skew-symmetric and therefore is sent to $\delta^{\mathrm{Ham}_{\pair{\cdot}{\cdot}},\nabla\!_W}_k$.

We refer to \cite{merval} for a detailed exposition on deformation complexes.\\
\end{remark}

\appendix
\section{The Infinitesimal Action and Divergences}\label{sec:infact}
This section is devoted to the further study of the interaction between the map $c_\nabla$ and a derivation action. Let $B$ be a $\mathbb{K}$-algebra and $M$ a left $B$-module. Recall the following in \cite{toyo}, Section 7:

\begin{definition}
A \textit{homological connection} on $M$ is a triple $\nabla = (P_\bullet,\partial_\bullet,\nabla_\bullet)$ where $P = (P_\bullet, \partial_\bullet)$ is a projective resolution of $M$ and $\nabla_i$ is a connection on $P_i$. A homological connection $\nabla$ is said to be flat if each $\nabla_i$ is flat.
\end{definition}

Now, we introduce the infinitesimal action.

\begin{definition}
Let $(\mathfrak d,\varphi,\rho)$ be a derivation action on $M$, and $\nabla$ as above. For $f\in\mathfrak d$, choose a lift $\lambda[f]$ of $\rho(f)\colon M \to M$ to the projective resolution $P$ once and for all, so that $\lambda[f]$ is a $\varphi(f)$-derivation and the association $f\mapsto \lambda[f]$ is $\mathbb{K}$-linear. We define the associated \textit{adjoint action} by
\[
	\ad_f\!\nabla = \big\{(L_{\varphi(f)}\otimes \id + \id\otimes \lambda[f]_i)\circ\nabla_i - \nabla_i\circ \lambda[f]_i \colon P_i \to \Omega^1B\otimes_BP_i\big\}_{i\geq 0}\,.
\]
Each map is an element of $\Omega^1(B,\End P_i)$.
\end{definition}

\begin{proposition}
$\Tr(\ad_f\!\nabla) \in |B|$ is independent of the choice of a lift.
\end{proposition}
\noindent Proof. If we take another lift $\lambda'[f]$, the difference is $B$-linearly null-homotopic: $\lambda[f] - \lambda'[f] = [\partial, h]$ for some $h\in\End_B^{(1)}(P)$. Then, we have
\begin{align*}
	\Tr&\{  \big(L_{\varphi(f)}\otimes \id + \id\otimes \lambda[f])\circ\nabla - \nabla\circ \lambda[f]\big) - \big(L_{\varphi(f)}\otimes \id + \id\otimes \lambda'[f])\circ\nabla - \nabla\circ \lambda'[f]\big) \}\\
	&= \Tr([[\partial, h],\nabla])\\
	&= \Tr([[\partial,\nabla],h] + [\partial, [h,\nabla]]).
\end{align*}
Since $[\partial,\nabla]$ and $[h,\nabla]$ are $B$-linear, $\Tr([[\partial,\nabla],h] + [\partial, [h,\nabla]])$ vanishes due to the cyclic invariance of the trace. (Note that $\nabla$ is not $B$-linear, so we cannot deduce $\Tr([[\partial, h],\nabla]) = 0$ immediately.)\qed\\

Now, we see the relation between the adjoint action and the divergence.

\begin{proposition}
Let $R = \nabla^2\in\Omega^2(M,\End P)$ be the curvature. Then, we have $\ad_f\!\nabla = D_\nabla(c_\nabla(f)) + i_{\varphi(f)}R$ in $\Omega^2(B,\End P)$ for $f\in\mathfrak d$.
\end{proposition}
\noindent Proof. We require to prove the equality $$(L_{\varphi(f)}\otimes \id + \id\otimes \lambda[f])\circ\nabla - \nabla\circ \lambda[f] = [\nabla, i_{\varphi(f)}\nabla - \lambda[f]] + i_{\varphi(f)}R.$$ The terms involving $\lambda[f]$ can be cancelled out, and the remaining is $(L_{\varphi(f)}\otimes \id)\circ\nabla = [\nabla, i_{\varphi(f)}\nabla] + i_{\varphi(f)}R$. 
Denoting $\varphi(f)$ by $f$, we have, for $p\in P$,
\begin{align*}
	\big((&L_f\otimes \id)\circ\nabla - [\nabla, i_f\nabla]\big)(p)\\
	&= L_f(\nabla'(p))\otimes \nabla''(p) - \nabla(i_f\nabla(p)) + (\id\otimes\,i_f\nabla)(\nabla(p))\\
	&= L_f(\nabla'(p))\otimes \nabla''(p) - \nabla(i_f\nabla'(p)\cdot \nabla''(p)) + (\id\otimes\,i_f\nabla)(\nabla'(p)\otimes \nabla''(p))\\
	&= L_f(\nabla'(p))\otimes \nabla''(p) - di_f(\nabla'(p))\cdot \nabla''(p) - i_f\nabla'(p)\cdot \nabla(\nabla''(p)) + \nabla'(p)\cdot i_f\nabla(\nabla''(p))\\
	&= i_fd(\nabla'(p))\cdot \nabla''(p) -  i_f(\nabla'(p)\cdot \nabla(\nabla''(p)))\\
	&= i_f(\nabla^2(p))\\
	&= i_f(R(p)).
\end{align*}
This completes the proof. \qed\\[-7pt]

\begin{corollary}
If $\nabla$ is trace-flat, i.e., $\Tr(R) = 0$, we have $\Tr(\ad_f\!\nabla) = d(\Div^{(\nabla,\varphi,\rho)}(f))$ in $\DR^1\!B$ for $f\in\mathfrak d$. 
\end{corollary}
\noindent Proof. By the proposition above, we have
\[
	\Tr(\ad_f\!\nabla) = \Tr(D_\nabla(c_\nabla(f))) + \Tr(i_{\varphi(f)}R).
\]
The second term on the right-hand side vanishes since the contraction and the trace map commute, and $\nabla$ is trace-flat. On the other hand, by Theorem \ref{thm:twdr}, we have $\Tr(D_\nabla(c_\nabla(f))) = d\Tr(c_\nabla(f)))$, which is further equal to $d(\Div^{(\nabla,\varphi,\rho)}(f))$ (up to a sign) by the definition of the divergence map in Section 7 of \cite{toyo}. This completes the proof. \qed\\[-7pt]

\begin{remark}
For an oriented Riemannian manifold $X$, the divergence $\sdiv_\mu(\xi)$ of a vector field $\xi$ on $X$ is a change of a fixed volume form $\mu$ by the Lie derivative: $\sdiv_\mu(\xi) \mu = L_\xi(\mu)$. This is recovered from our definition of the divergence by setting $\mathbb{K} = \mathbb{R}$, $B = C^\infty(X)$, $M = \wedge^\mathrm{top}T^*X$, $\mathfrak d = \Der(C^\infty(X))$,
\begin{align*}
	\varphi = \id&\colon \mathfrak d \to \Der_\mathbb{K}(B),\\
	\rho = L&\colon \mathfrak d \to \End_\mathbb{K}(M),\textrm{ and }\\
	\nabla=\nabla\!_\mu&\colon M \to \Omega^1B\otimes_BM\colon \mu \mapsto 0.
\end{align*}
Here, the flat connection $\nabla\!_\mu$ should be regarded as an incarnation of the volume form $\mu$. The Proposition above says that, under the assumption of the flatness, the infinitesimal change of the connection (conversely thought of as a non-commutative volume form) by a vector field $f$ roughly gives the divergence as one may expect.
\end{remark}

The map $c_\nabla$ also naturally emerges from the infinitesimal adjoint action.
\begin{proposition}
For any homological connection $\nabla$ and $f\in \mathfrak d$, $\iota(\ad_f\!\nabla)$ is equal to
\[
	\ad_f(c_\nabla) := [\lambda[f],c_\nabla] - c_\nabla\circ\ad_f
\]
in $C^1_\mathrm{CE}(\mathfrak d, \End_B(P))$ up to homotopy term.
\end{proposition}
\noindent Proof. For $g\in\mathfrak d$, we have
\begin{align*}
	 \iota(\ad_f\!\nabla)(g) &= i_g\big((L_f\otimes \id + \id\otimes \lambda[f])\circ\nabla - \nabla\circ \lambda[f])\big)\\
	 &= (i_gL_f\otimes\id + \,i_g\otimes\lambda[f])\circ\nabla - i_g\nabla\circ\lambda[f].
\end{align*}
On the other hand, we have
\begin{align*}
	([\lambda[f],c_\nabla] - c_\nabla\circ\ad_f)(g) &= [\lambda[f] ,i_g\nabla - \lambda[g]] - (i_{[f,g]}\nabla - \lambda[[f,g]])\\
	&=  [\lambda[f] ,i_g\nabla] - i_{[f,g]}\nabla + ( \lambda[[f,g]] -  [\lambda[f] ,\lambda[g]] ).
\end{align*}
The last two terms $\lambda[[f,g]] -  [\lambda[f] , \lambda[g]]$ is equal to $[\partial,h]$ for some $h\colon \mathfrak d \to \End^{(-1)}_B(P)$ since it is the difference of the lifts of the same map $\rho([f,g]) = [\rho(f),\rho(g)]$. Therefore,
\begin{align*}
	 (\iota(\ad_f\!\nabla) - \ad_f(c_\nabla))(g) &= (i_gL_f\otimes\id + \,i_g\otimes\lambda[f])\circ\nabla - i_g\nabla\circ\lambda[f] - [\lambda[f] ,i_g\nabla] + i_{[f,g]}\nabla - [\partial,h]\\
	 &= (L_fi_g\otimes\id +\,  i_g\otimes\lambda[f])\circ\nabla  - \lambda[f]\circ i_g\nabla - [\partial,h]\\
	 &= - [\partial,h].
\end{align*}
The last equality comes from the assmuption that $\lambda[f]$ is a $\varphi(f)$-derivation.\qed\\

\small
\bibliographystyle{alphaurl}
\bibliography{faoetc.bib}

\end{document}